\documentclass[12pt]{article}
\usepackage{latexsym,amssymb,amsmath,amsthm,amsxtra,amsfonts}
\usepackage{graphicx}
\usepackage{enumerate}
\newtheorem{thm}{Theorem}[section]
\newtheorem{prop}[thm]{Proposition}
\newtheorem{lem}[thm]{Lemma}

\newtheorem*{nnthm}{Theorem}
\newtheorem{conj}{Conjecture}
\theoremstyle{definition}
\newtheorem{defn}[thm]{Definition}

\theoremstyle{remark}
\newtheorem{rmk}{Remark}

\newtheorem*{rmk*}{Remark}
\newtheorem*{rmks*}{Remarks}
\newenvironment{enumeratei}{\begin{enumerate}[\upshape (i)]}
                           {\end{enumerate}}
\newcommand{\deq}{\overset{{\rm def}}{=}}            
\newcommand{\E}{{\mathbb E}}                
\newcommand{\Eb}[1]{\E\left[ #1 \right]}
\renewcommand{\P}{{\mathbb P}}              
\newcommand{\Pb}[1]{\P\left[ #1 \right]}
\newcommand{\N}{{\mathbb N}}

\newcommand{\Z}{{\mathbb Z}}
\newcommand{\Q}{{\mathbb Q}}

\newcommand{\R}{{\mathbb R}}

\newcommand{\given}{\mid}

\newcommand{\Lp}[2]{\|#2\|_{#1}}

\renewcommand{\and}{\; \mbox{ and } \;}
\newcommand{\iid}{i.i.d.\ }

\DeclareMathOperator{\cov}{Cov}

\newcommand{\sequ}[1]{\{ #1_n \}_{n=0}^\infty}
\newcommand{\fc}{\frac{1}{2\pi}}

\newcommand{\ep}{\epsilon}
\newcommand{\al}{\alpha}

\renewcommand{\phi}{\varphi}
\newcommand{\Qs}{{\mathcal Q}}

\newcommand{\bias}{\theta}

\newcommand{\rhs}{right-hand side }
\newcommand{\prm}{{\mathcal P}}

\newcommand{\chf}{\psi}
\newcommand{\rw}{S}
\newcommand{\bo}[1]{{\boldsymbol #1}}
\newcommand{\av}[2]{{\overline{#1}({#2})}}
\newcommand{\lbr}{[}
\newcommand{\rbr}{]}

\newcommand{\vect}[1]{{\boldsymbol #1}}
\newcommand{\tens}[1]{{\boldsymbol #1}}

\newcommand{\M}{\tens{M}}

\newcommand{\inv}[1]{#1'}
\begin{document}
\begin{center}
{\sc \large Identifying  several biased coins  \\
encountered by a hidden random walk }
%
\vspace{.25in}

{\sc  
David A.~Levin\footnote{Department of Mathematics, University of 
Utah, Salt Lake City, UT 84112. {\tt levin@math.utah.edu}.} 
and Yuval Peres\footnote{Departments of Statistics and Mathematics, 
University of California, Berkeley, CA 94720.
{\tt peres@stat.berkeley.edu}.
{\tt http://stat-www.berkeley.edu/}$\tilde{\ }${\tt peres}.
Research supported in part by NSF Grants
DMS-0104073 and DMS-0244479.
}}

\vspace{0.25in}
\sc \today \rm
\end{center}
\begin{abstract}
Suppose that attached to each site $z \in \Z$ is a coin
    with bias $\bias(z)$, and only finitely many of these
    coins have non-zero bias.  Allow a simple random walker to generate
    observations by tossing, at each move, the coin attached
    to its current position.  Then we can determine
    the biases $\{\theta(z)\}_{z \in \Z}$, {\em using only the
    outcomes of these coin tosses and no information about the
    path of the random walker}, up to a shift and reflection of $\Z$.  
    This generalizes a result of Harris and Keane.
\end{abstract}



\section{Introduction}

A {\em coin toss} with {\em bias} $\theta$ is a $\{-1,1\}$-valued
random variable with expected value $\theta$, and a {\em fair coin} is a coin
toss with bias zero. 

Harris and Keane \cite{HK:RCT} introduced a model for
sequences of coin tosses with a particular kind of dependence.
Let $\rw = \{\rw_n\}_{n=0}^{\infty}$ be simple random walk on $\Z$.
Suppose that whenever $\rw_n=0$,
 an independent coin with bias $\bias \ge 0$ is tossed,
while at all other times an independent fair coin is tossed. 
Suppose that we are given $X=\{X_n\}_{n=0}^{\infty}$,
 the record of coin tosses obtained, but the path of the walker,
 $\rw$, is hidden from us.
Can $X$ be distinguished from 
a sequence of \iid fair coin tosses? \newline
If it can be distinguished, can the parameter $\bias$ be recovered, again
only from $X$?

Harris and Keane proved the following:

\begin{nnthm}[{\cite[Theorem 2]{HK:RCT}}] 
There exist a sequence of functions $\Theta_n: \{-1,1\}^n \rightarrow \R$, 
not depending on $\theta$, so that 
$\lim_{n \rightarrow \infty} \Theta_n(X_1,\ldots,X_n) = \bias$ almost surely.
\end{nnthm}

In fact, the model in \cite{HK:RCT} is more general, allowing
the walk to be  any null-recurrent Markov chain.  Let $\mu_{\bias}$ be
the distribution of $X$ on $\{-1,1\}^{\N}$. Recall that two Borel measures 
are {\em mutually singular} if there is a Borel set which is null for
one measure such that its complement is null for the other measure;
in this case the measures are distinguished 
by almost any observation.  Whether or not
$\mu_{\bias}$ and $\mu_{0}$ are mutually  singular
depends on the graph, and as was shown in
\cite{LPP:PTRCT}, sometimes on the bias $\bias$. 
\cite{HK:RCT} provides a partial criterion for determining 
singularity based on the graph, while the full story -- in particular
the role of $\bias$ for certain graphs -- is completed in
\cite{LPP:PTRCT}.

In this paper, this model is generalized in a different
direction.  Label each $z \in \Z$ with $\bias(z) \in [0,1]$. 
Now allow a random walker to move on $\Z$,
and suppose that on her $n$th move she is at position $S_n$ and 
she tosses a coin with bias
$\bias(S_n)$, generating the random variable $X_n$. 
Again, the path of the walker, $\{S_n\}_{n=0}^\infty$ is not observed,
and the output is only the outcomes of the coin tosses,
that is, $\{X_n\}_{n=0}^\infty$.  The main result of this
paper says that the biases $\{\bias(z) : z \in \Z \}$ can
be recovered, up to a shift and reflection of $\Z$, in the case when the
number of vertices with $\bias(z) \neq 0$ is finite:

\begin{thm} \label{thm:coins}
    Suppose that attached to each site $z \in \Z$ is a coin
    with bias $\bias(z)$, and only finitely many of these
    coins have non-zero bias.  Allow a simple random walker to generate
    observations by tossing, at each move, the coin attached
    to its current position.  Then it is possible to determine
    the biases $\{\theta(z)\}_{z \in \Z}$, 
  up to a shift and reflection of $\Z$,  {\em using only 
    a single observation of the infinite sequence of
    these coin tosses and no information about the
    path of the random walker}.

\end{thm}

\begin{figure}
    \begin{center}
        \includegraphics{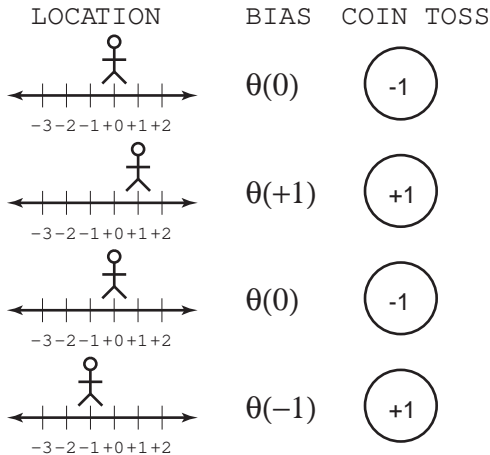}
    \end{center}
    \caption{A walker performing his coin flips.  We only see the
    right-hand column, but would like to recover the center
    column!}
\end{figure}

Another related model was introduced by Benjamini and Kesten \cite{BK:Scen}
(and independently, by Den Hollander and Keane).
Instead of attaching a coin with bias $\theta(z)$ to
each $z \in \Z^d$, a fixed {\em color} is assigned to each
vertex.  The random walker reports only the sequence of
colors visited, while just as in our generalization
of the Harris-Keane model, the path itself
is hidden.  
Further
investigations of the model introduced in \cite{BK:Scen} can be found
in  \cite{H:PS,H:RWS,L:IS,LM:2d,LM:Corr,Ma,MR:Poly,MR:Err} and 
a survey of results (up to 1998) can be found in \cite{K:Scenery}.

This paper is organized as follows.  Theorem \ref{thm:coins}
holds more generally, where the observed random variables 
are not restricted to be coin tosses.  Some definitions and
the statement of the more general Theorem \ref{thm:main}
can be found in Section \ref{sec:defns}.  The main requirements of 
the proof of this theorem are a strong law of large numbers, 
and some algebraic results, found in Section \ref{sec:propositions}. 
Because the special case of coin tosses is more readily understood,
a separate proof of Theorem \ref{thm:coins} is supplied,
after which Theorem \ref{thm:main} is proven, both
in Section \ref{sec:proof}.  Finally, some unsolved problems
are mentioned in Section \ref{sec:conclusion}.

\section{Definitions and Statement of Result} \label{sec:defns}
Let $q$ be a probability mass function on $\Z$, and
$q^{*n}$ its $n$-fold convolution.
Assume that
\begin{equation} \label{eq:qconds1}
\sum_{z \in \Z} z q(z) = 0 \quad \text{and} \quad
\sum_{z \in \Z} z^2 q(z) < \infty.
\end{equation}
Also, we will assume that $q$ is symmetric:
\begin{equation} \label{eq:symmetric}
    q(z) = q(-z)\, \text{ for all } z \in \Z \,.
\end{equation}
In addition, to simplify the exposition we will
suppose that
\begin{equation} \label{eq:qcondsnp}
{\rm g.c.d.}\{n : q^{*n}(0) > 0  \} = 1 \,.
\end{equation}
\begin{rmk*}
Condition (\ref{eq:qcondsnp}) rules out, for example, 
$q = \frac{1}{2}\delta_{-1} + \frac{1}{2}\delta_{+1}$,
the mass function corresponding to the increments of
simple random walk on $\Z$.  However, the results in this
paper (in particular, Theorem \ref{thm:main}) can be
obtained for periodic random walks by considering 
the process observed at times $\{rn+k : n \in \N\}$,
where $r$ is the period (the left-hand side of \eqref{eq:qcondsnp}).
\end{rmk*}

The random walk with increments distributed 
according to $q$ will be denoted $\rw = \sequ{\rw}$:
\begin{equation*}
\P[S_{n+1} - S_{n} = z \given
S_0, \ldots, S_n ]  =  q(z) \,.
\end{equation*}
We use the notation 
$P^n(z,w) \deq  \P[\rw_{n} = w \given \rw_0 = z]$
and $u_n  \deq  P^n(0,0)$.
\begin{defn}
    \rm A {\em stochastic scenery} is a 
    map $\eta:\Z \rightarrow \prm(\R)$, where
    $\prm(\R)$ is the set of all probability
    measures on $\R$.  
\end{defn}
We assume that
there is a known reference
distribution $\alpha \in \prm(\R)$ which appears everywhere but finitely many
vertices: that is,
\begin{equation} \label{eq:finassum}
|\eta^{-1}(\prm \setminus \{ \alpha \})| < \infty\,. 
\end{equation}
\begin{defn}
    The {\em sequence of observations induced by $\rw$
    and $\eta$} is the stochastic
    process $X = \sequ{X}$, whose law conditioned
    on $\rw$ has independent marginal distributions given by
    $\P [X_n \in \cdot \given \rw] = \eta(\rw_n)$.
\end{defn}

\begin{thm} \label{thm:main} 
For a random walk $\rw$ with symmetric
increment distribution $q$ satisfying
(\ref{eq:qconds1}) and (\ref{eq:qcondsnp}),
and stochastic scenery $\eta$ obeying (\ref{eq:finassum}), 
let $X = \sequ{X}$ be the 
observations induced by $\rw$ and $\eta$.
Then using only the observations $\{X_{n}\}_{n=0}^\infty$, we can determine
a scenery $\hat{\eta}$ so that
$$
\hat{\eta} = \eta \circ T \; \text{a.s.}\,,
$$
where $T:\Z \rightarrow \Z$ is 
either a shift or a shift composed with a reflection,
and may depend on $\eta$.
\end{thm}
\begin{rmk*}
As noted by a referee, if we are given an infinite sequence of
observations then we do not need to be given the reference
distribution $\alpha$, as it can be a.s.\ reconstructed from the
observations.
\end{rmk*}
\section{A few ingredients} \label{sec:propositions}
%
%
In this section, we first state Proposition \ref{prop:sl},
which shows 
how the observations $\{X_{n}\}$ can
be used to determine a family of parameters $\{p_{\vect{t}}\}$.
We show in Proposition \ref{prop:ptoq}
that we can express a new family of parameters $\{Q^k_{\vect{d}}\}$
as functions of $\{p_{\vect{t}}\}$.  These new parameters 
determines the scenery $\eta$, and hence the scenery can be
recovered from the observations $\{X_n\}$. 
(This is the substance of 
the proofs of Theorem \ref{thm:coins} and Theorem \ref{thm:main}
in Section \ref{sec:proof}.)
We conclude the present section with a proof of Proposition
\ref{prop:sl}.

\begin{defn} \label{def:p}
\begin{enumeratei}
\item
    For a vector $\vect{t} \in \N^{k}$
    (where $k \in \N$),
    and a bounded and
    measurable $\phi:\R \rightarrow \R$,
    define
    \begin{equation*}
    p_{\vect{t}}(\phi)
     \deq  \sum_{z\in \Z} \E [ 
    \phi(X_0)\phi(X_{t_1}) \phi(X_{t_1+t_2}) \cdots 
        \phi(X_{t_1 + \cdots + t_k}) \given \rw_0 = z ].
    \end{equation*}
\item
    For $r \in \N$, let $\vect{p}^k(\phi)$ be
    the vector with components 
    $$\{p_{\vect{t}}(\phi) \; : \;
    \vect{t} \in \{1,\ldots,r\}^k\ \} \,.
    $$
\end{enumeratei}
\end{defn}
For measurable functions $\phi$, we will use the supremum norm
$\| \phi \|_{\infty}$.  When we know that all the measures $\eta(z)$
are supported on a common set $E$, we will take this norm to be the
essential supremum over $E$.  For example, when we are dealing
with coin tosses only, $\| \phi \|_\infty \leq 1$ means
that $\phi$ is bounded by $1$ on $\{-1,1\}$, a condition satisfied by
$\phi(x)=x$.
\begin{prop} \label{prop:sl}
    Fix a measurable function $\phi$ with 
    $\| \phi \|_{\infty} \leq  1$ 
    and
    $\int_{\R} \phi d \alpha = 0$,
    and a vector $\vect{t} \in \N^{s}$.
    Then there exists a sequence $\{g_N\}$ of measurable functions,
    $
    g_N : \R^N \rightarrow \R \,,
    $
    not depending on the stochastic scenery $\eta$, with 
\begin{equation*}
\lim_{N \rightarrow \infty} g_N(X_1,\ldots,X_N) =
    p_{\vect{t}}(\phi)
    \quad \text{almost surely}. 
\end{equation*}
\end{prop}

The information in the collection $\{p_{\vect{t}}\}$
must be transformed into a more useful form for
us to recover the stochastic scenery $\eta$.
We provide now several results enabling this
transformation.

\begin{lem} \label{lem:fullrank}
    For a given transition semigroup $P^t$, 
    and $m,r \in \N$, define the matrix $\M$ as
    \begin{equation} \label{eq:promatrix}
    \M  \deq \left[
    \begin{array}{llll}
        P^{1}(0,0) & P^{1}(0,1) & \cdots & P^{1}(0,m) \\
        P^{2}(0,0) & P^{2}(0,1) & \cdots & P^{2}(0,m) \\
        \vdots & \vdots & \ddots & \vdots \\
        P^{r}(0,0) & P^{r}(0,1) & \cdots & P^{r}(0,m)
    \end{array}
    \right] \,.
    \end{equation}
    Let $\rw$ be an aperiodic and symmetric random walk on $\Z$, whose 
    increments have probability mass function $q$: that is,
        $\P[\rw_{t+1} = z \given \rw_t = w] = q(z-w)$.
    \begin{enumeratei}
        \item 
        If 
        $P^t$ is the transition semigroup for
        $\rw$, then
        $\M$ has rank $(m+1)$ for some $r$.
        \item
        If $\rw$ is the {\em simple} random walk on $\Z$,
        and
        $P^{t}(0,z) = \P[\rw_{2t}=2z]$, 
        then any $(m+1) \times (m+1)$
        square submatrix of $\M$ is non-singular.
   \end{enumeratei}
\end{lem}
\begin{proof}
If $\chf$ is the characteristic function of $q$, so
that $\chf^t$ is the characteristic function of the
$t$th convolution power of $q$, then 
by Fourier inversion
\begin{equation}\begin{split}\label{eq:fourinv}
\fc \int_{-\pi}^\pi \cos(uj) \chf^t(u) du & = 
2\P[ \rw_t = j]. 
\end{split}\end{equation}
Suppose that for each $r$, the $m+1$ columns of $\M$ are 
linearly dependent, that is, there exist $\{a_j(r)\}_{j=0}^m$ 
(not all zero) so that
$$
\sum_{j=0}^m a_j(r) P^t(0,j) = 0 \text{ for all } t \leq r \,.
$$
By rescaling, it can be assumed that $\sum_{j=0}^m a_j(r)^2 = 1$.
Thus there exist subsequential limits $\{a_j\}_{j=0}^m$ (not all
zero) so that
\begin{equation} \label{eq:lindepend}
\sum_{j=0}^m a_j P^t(0,j) = 0 \text{ for all } t \,.
\end{equation}
Hence, the matrix $\M$ in (\ref{eq:promatrix}) fails to have rank $(m+1)$ for 
any $r$ if and 
only if there exist $\{a_j\}_{j=0}^{m}$ 
satisfying \eqref{eq:lindepend}.
By (\ref{eq:fourinv}), this is the same as
\begin{equation} \label{eq:lindep}
\int_{-\pi}^\pi g(u) \chf^t(u) du = 0
\text{ for all }  t,
\end{equation}
where $g(u) = \sum_{j=0}^{m} a_j \cos(uj)$.
We now show that (\ref{eq:lindep}) cannot hold:

Observe that $g$ is an even function that is non-zero in a neighborhood of $0$,
except possibly at $0$.
It can be assumed without loss of generality (by multiplying the
coefficients $a_j$ by $-1$) that $g$ is non-negative near $0$.
Also $\chf$ is even, has value $1$ at the origin, and
is positive in a neighborhood $[-\delta,\delta]$ of $0$.  

Choose $\delta_1 < \delta$ so that
\begin{equation*}
    B \deq \min_{|u| \leq \delta_1} \chf(u) > \max_{\delta < |u| \leq \pi}
    \chf(u) \deq A 
\end{equation*}
This is possible because, as $\chf$ is a characteristic function of an
aperiodic distribution, $A < 1$.
Choose $\delta_1$ even smaller to ensure $\min_{\frac{\delta_1}{2}\leq |u| \leq \delta_1}g(u) = \epsilon > 0$.
Hence
\begin{align} 
    \int_{-\pi}^\pi g(u) \chf^t(u) du
    & = \int_{-\delta}^{\delta} g(u) \chf^t(u) du
    + \int_{[-\pi,\pi]\setminus [-\delta,\delta]}  g(u) \chf^t(u) du 
    \nonumber \\
    & \geq \int_{\delta_1/2}^{\delta_1} g(u) \chf^t(u) du
    + \int_{[-\pi,\pi]\setminus [-\delta,\delta]}  g(u) \chf^t(u) du
    \nonumber \\
    & \geq  \frac{\delta_1}{2}\ep B^t - \| g \|_{\infty} 2\pi A^t \,.
    \label{eq:posgeo}
\end{align}

For $t$ large enough, the \rhs of (\ref{eq:posgeo})
is positive, contradicting (\ref{eq:lindep}).
Thus the matrix (\ref{eq:promatrix}) has rank $(m+1)$ when 
$r$ is sufficiently large, proving (i).

We now discuss the special case of simple random walk. 
Suppose there exist $a_0,\ldots,a_m$ and $t_0,t_1,\ldots,t_m \in 2\Z$ so that
for $t=t_0,\ldots,t_m$,
\begin{equation*}
    \sum_{k=0}^m a_k \P[S_{2t}=2k]  =  0.
\end{equation*}
Then 
\begin{equation} \label{eq:lincomb}
\sum_{k=0}^m a_k \frac{(2t)!}{(t-k)!(t+k)!} = 0.
\end{equation}
Multiplying both sides of (\ref{eq:lincomb}) by 
$\frac{(t+m)!t!}{(2t)!}$ yields
\begin{equation} \label{eq:npoly}
0 = \sum_{k=0}^m a_k (t+m)\cdots(t+k+1) t \cdots (t-k+1).
\end{equation}
If $f_k(x) \deq (x+m) \cdots (x+k+1) x \cdots (x-k+1)$, then  
(\ref{eq:npoly}) shows that 
\begin{equation} \label{eq:polyzeros}
f(x) \deq \sum_{k=0}^m a_k f_k(x) = 0 \; \mbox{ for } \; 
x = t_0,\ldots,t_m \,.
\end{equation}
Since $f$ is a polynomial of degree $m$, and by (\ref{eq:polyzeros}) 
$f$ has $m+1$ zeros,  
it follows that $f \equiv 0$.
But $f(-m)=0$ implies $a_m=0$. 

Multiplying both side of (\ref{eq:lincomb})
by $\frac{(t+m-1)!t!}{(2t)!}$ and arguing as above
shows that $a_{m-1}=0$.   Repeating this procedure shows
that $a_k = 0$ for all $k$, establishing (ii).
\end{proof}
\begin{rmk} \label{rmk:m_suppressed}
    We have suppressed the dependence of the $r \times (m+1)$ matrix $\M$ on
    both $m$ and $r$.  For each $m$, fix an $r_m$ so that
    $\M$ has rank $(m+1)$.  In the sequel, the matrix $\M$
    will depend only on $m$, having dimension $r_m \times 
    (m+1)$.  We will write $r$ for $r_m$.  
\end{rmk}
We now introduce some notation:
\begin{defn} \label{defn:func_ave}
    For any function $\phi:\R \rightarrow \R$,
    define
    \begin{equation} \label{eq:func_ave}
        \av{\phi}{z} \deq \int_\R \phi d \eta(z) \,,
        \quad \mbox{and} \quad
        \phi_0 \deq \phi - \int_\R \phi d \alpha \,.
    \end{equation}
    (Recall that $\eta(z) \in \prm(\R)$.)
\end{defn}
The parameters $\{\vect{Q}^{k}(\phi)\}$ which will be used in
the reconstruction of the scenery are now defined:
\begin{defn}    \label{defn:q}
    For $\phi$ a bounded measurable function and $k \in \N$, define the
    vector $\vect{Q^k(\phi)} \in \R^{[m]^k}$, indexed
    by $\vect{d} \in [m]^{k}$, by
    \begin{equation} \label{eq:qdefn}
       Q^{k}_{\vect{d}}(\phi) \deq
       \sum_{\substack{(z_1,\ldots,z_{k+1}) \\ |z_j-z_{j-1}|=d_j}}
       \av{\phi}{z_1}\av{\phi}{z_2}\cdots\av{\phi}{z_{k+1}}.
   \end{equation}
(Here $[m] \deq \{0,1\ldots,m\}$.)
\end{defn}

Assumption \eqref{eq:finassum} of Theorem \ref{thm:main}
requires that the stochastic
scenery $\eta$ differs from a reference distribution
$\alpha$ at only finitely many vertices.  With this
in mind, we make the following definitions:
\begin{defn} \label{defn:abl}
Let
\begin{equation} \label{eq:abdef}
a \deq \inf\{ j : \eta(j) \neq \alpha\}, \quad  
b \deq \sup\{ j : \eta(j) \neq \alpha\}, \quad
\text{and} \quad \ell \deq b - a.
\end{equation}
\end{defn}
\noindent By assumption (\ref{eq:finassum}),
$- \infty < a \leq b < \infty$.

Finally, we state the main algebraic result used in the proof of
Theorem \ref{thm:main}:
\begin{prop} \label{prop:ptoq} 
If $m \geq \ell$ 
and $\int \phi d \alpha = 0$,
then $\vect{Q^{k}}(\phi) \in \R^{[m]^k}$ 
is a linear transformation 
of $\vect{p}^k(\phi) \in \R^{r^k}$ 
(defined in Definition \ref{def:p}.)
\end{prop}
\proof
Let us begin with the case of $k=1$.
Let $\M$ be the matrix in \eqref{eq:promatrix}, and
$\inv{\M}$ its left inverse, shown to exist in Lemma
\ref{lem:fullrank}.

We have
\begin{align}  
p_{t}(\phi)
& = \sum_{z \in \Z} \E [ \phi(X_t) \phi(X_0) \given \rw_0 = z] 
\nonumber \\
& = \sum_{z \in \Z} \sum_{w \in \Z} \E [ \phi(X_t) \phi(X_0) \given \rw_0 = z,
 \rw_t = w] \P [ \rw_t = w \given \rw_0 = z] \nonumber \\
& = \sum_{d \in \Z} P^t(0,d) \sum_{\substack{(z,w)\\ z-w=d}} 
\av{\phi}{z}\av{\phi}{w} \,.\label{eq:polyformgen}
\end{align}
Because $q$ is symmetric, we have $P^t(0,d) = P^t(0,-d)$, and hence
{}from (\ref{eq:polyformgen}) we conclude that
\begin{equation*}  
p_{t}(\phi)
  = \sum_{d=0}^{\infty} P^t(0,d) \sum_{\substack{(z,w) \\ |z-w|=d}}
\av{\phi}{z}\av{\phi}{w} \,.
\end{equation*}
If $d \geq m$, then $d \geq \ell$ and the factor $\sum_{|z-w|=d}
\av{\phi}{z}\av{\phi}{w}$ vanishes, so we have
\begin{equation} \label{eq:tosolve}
    p_{t}(\phi) = \sum_{d=0}^{m} M_{t,d} \cdot 
    Q^1_d(\phi) \,.
\end{equation}
Using Definition \ref{def:p}, we rewrite the collection of
$r$ equations \eqref{eq:tosolve} as
\begin{equation} \label{eq:tosolvevec}
    \vect{p}^1(\phi) =
    \M \vect{Q^1(\phi)}.
\end{equation}
Multiplying both sides of \eqref{eq:tosolvevec} on the
left by $\inv{\M}$ yields
$\vect{Q^1(\phi)} = \inv{\M}\vect{p}^1(\phi)$.

Now consider the case $k=2$.  We have
\begin{align}  
p_{(t_1,t_2)}(\phi)
& =  \sum_{z,w,v \in \Z} 
\E [ \phi(X_{t_1 + t_2})\phi(X_{t_1})\phi(X_{0}) \given
S_0  = z, S_{t_1} = w, S_{t_2} = v]  \nonumber \\
& \quad \times P^{t_1}(0,w-z)P^{t_2}(0,v-w) \nonumber \\
& =  \sum_{z,w,v \in \Z} \av{\phi}{z}\av{\phi}{w}\av{\phi}{v}
P^{t_1}(0,w-z)P^{t_2}(0,v-w)  \nonumber \\
& =  \sum_{d_1=-\infty}^\infty \sum_{d_2=-\infty}^\infty
P^{t_1}(0,d_1)P^{t_2}(0,d_2) \sum_{
\substack{(z,w,v)\\ w-z = d_1, v-w=d_2}} \av{\phi}{z}\av{\phi}{w}
\av{\phi}{v} \nonumber \\
& =  \sum_{d_1=0}^\infty \sum_{d_2=0}^\infty 
P^{t_1}(0,d_1)P^{t_2}(0,d_2) 
\sum_{\substack{(z,w,v) \\ |w-z| = d_1, |v-w| = d_2}}
\av{\phi}{z}\av{\phi}{w}\av{\phi}{v}. \nonumber
\end{align} 
Again, since $\sum_{|z-w| = d_1, |w-v| = d_2}
\av{\phi}{z}\av{\phi}{w}\av{\phi}{v}$ vanishes for
$d_1 \vee d_2 > m \geq \ell$, we have
\begin{equation} \label{eq:polythrees}
p_{(t_1,t_2)}(\phi) =   \sum_{d_1=0}^m \sum_{d_2=0}^m 
P^{t_1}(0,d_1)P^{t_2}(0,d_2) 
\sum_{\substack{(z,w,v) \\ |z-w| = d_1, \\|w-v| = d_2}}
\av{\phi}{z}\av{\phi}{w}\av{\phi}{v}  \,.
\end{equation}
Recall that the tensor product 
of the matrix $\M$ 
in \eqref{eq:promatrix} with itself is given by 
\begin{equation*}
\M \otimes \M =
\left[
\begin{array}{lcl}
M_{1,0} \M & \cdots & M_{1,m} \M \\
\vdots & \ddots & \vdots \\
M_{r,0}\M & \cdots & M_{r,m} \M
\end{array}
\right].
\end{equation*}
There are $r^2$ rows, indexed by 
$\vect{t} \in \{(t_1,t_2) \; : \;
1 \leq t_1, t_2 \leq r\}$, and $(m + 1)^2$ columns indexed
by $\vect{d} \in [m]^2 = \{(d_1,d_2) \; : \; d_i = 0,1,\ldots,m \}$.
Each row is of the form
\begin{equation*}
    \left(P^{t_1}(0,d_1)P^{t_2}(0,d_2)\right)_{d_1,d_2=0}^{m}. 
\end{equation*}
Thus we can rewrite the equations \eqref{eq:polythrees} for
$\vect{t} \in [r]^2$ compactly as
\begin{equation} 
    \vect{p}^2(\phi) = \left(\M \otimes \M \right) 
    \vect{Q^2(\phi)} \,. \nonumber 
\end{equation}
Because each $\M$ has rank $(m+1)$, the matrix
$\M \otimes \M$ has rank $(m+1)^{2}$ (see
\cite[Theorem 4.2.15]{HJ:TMA}). Thus there is a left inverse 
$\inv{(\M \otimes \M)}$ to $\M \otimes \M$ and we can write
\begin{equation} 
    \vect{Q^2(\phi)} = \inv{(\M \otimes \M)} \vect{p}^2(\phi)
    \,. \nonumber 
\end{equation}
In general, $\M^{k \otimes}$ has rank $(m+1)^{k}$, and hence has
a left inverse, enabling us to write
\begin{equation} \label{eq:qfromp}
    \vect{Q^k(\phi)} = \inv{\left({\M}^{k \otimes}\right)} \vect{p}^k(\phi).
\end{equation}
\endproof
%
%
\begin{rmk}
    The vectors $\vect{p}^k(\phi)$ depend on the constant
    $r$ (see Definition \ref{def:p}), which we have
    taken to be the constant $r_m$ defined in 
    Remark \ref{rmk:m_suppressed}.  Thus these
    vectors depend on $m$, although we have
    again suppressed this dependence in the notation.
\end{rmk}
\begin{lem} \label{lem:balem}
    Let $\ell$ be the constant given in Definition
    \ref{defn:abl}.  
    Then there exists a sequence of random variables
    $\{\ell_n\}_{n=0}^\infty$, not depending on $\eta$, 
    with $\ell_n$ measurable with 
    respect to $\sigma(X_1,\ldots,X_n)$, and so that 
    $\ell_n=\ell$ eventually, a.s.
\end{lem}
Recall that a class $\Phi$ of measurable functions on $\R$ is said to be
{\em measure determining} if
$\int_{\R} \phi d\mu  = \int_{\R} \phi d\nu$
for each $\phi \in \Phi$ implies $\mu=\nu$.
Let $\Phi$ be a countable
measure determining class of functions on $\R$
with $\Lp{\infty}{\phi} \leq 1$.
We will enumerate $\Phi$ as $(\phi_1,\phi_2,\ldots)$.
Also define
$\Phi_1 \deq \left\{ \frac{f + cg}{\|f+cg\|_{\infty}} \; : \; f,g \in
\Phi,\, c \in \Q \right\}$, 
and enumerate its elements as $(h^1,h^2,\ldots)$.

\begin{lem} \label{lem:disting}
    Let $\mu$ and $\nu$ be two probability measures, neither
    equal to $\alpha$. Recall that $\psi_0 = \psi - \int \psi d\alpha$.
        If $\mu \neq \nu$, there exists $\psi \in \Phi_1$
        with 
        \begin{equation} \label{eq:clpsi}
            \int \psi_0 d\mu \neq 0\,, \quad
            \int \psi_0 d\nu \neq 0\,, \quad \mbox{and} \quad
            \int \psi_0 d\mu \neq \int \psi_0 d\nu \,.
        \end{equation}
\end{lem}
\proof[Proof]

Because $\mu \neq \nu$ and $\Phi$ is measure determining,
there is an 
$f \in \Phi$ so that
$\int f d\mu \neq \int f d\nu $.
By subtracting $\int f d\alpha$, we have
$\int f_0 d\mu \neq \int f_0 d\nu $.
Then 
one of $\int f_0 d\mu$ and $\int f_0 d\nu $ is
nonzero; without loss of generality assume
$\int f_0 d\mu \neq 0$.
Similarly, there
is a $g \in \Phi$ so that $\int g d\nu \neq \int g d\alpha$, that is,
$\int g_0 d\nu \neq 0$.

Hence, for small enough $c \in \Q$, we have

\begin{equation} \label{eq:linnz}
    \int f_0 d\mu + c \int g_0 d\mu \neq 0 \,,
    \text{ and }
    \int f_0 d\nu + c \int g_0 d\nu \neq 0 \,.
    \end{equation}
By taking $c$ even smaller, since $\int f_0 d\mu \neq \int f_0 d\nu$,
we have
\begin{equation} \label{eq:lindif}
    \int f_0 d \mu + c \int g_0 d \mu \neq \int f_0 d\nu + c\int g_0 
    d\nu \,.
    \end{equation}
Thus, if $\psi = \frac{f + cg}{\|f +cg\|_{\infty}}$, so that
$\psi_0 = \frac{f_0 + cg_0}{\|f + cg\|_{\infty}}$, 
we have from (\ref{eq:linnz}) and
(\ref{eq:lindif}) that $\psi$ obeys (\ref{eq:clpsi}).
\endproof

We can now proceed with the proof of Lemma \ref{lem:balem}:
\begin{proof}[Proof of Lemma \ref{lem:balem}]
As before, 
let $\inv{\M}$ be the left inverse to the matrix $\M$
in \eqref{eq:promatrix}.  

Let $\ell^h \deq \sup\{|z-w| \;:\;
\av{h_0}{z}\av{h_0}{w} \neq 0\}$.  Recall that $(h^1,h^2,\ldots)$ is an
enumeration of $\Phi_1$, and write $\ell^j$ for $\ell^{h^j}$.
Lemma \ref{lem:disting} guarantees that $\sup_{j \leq r} \ell^j = \ell$ for 
$r$ sufficiently large.  
Define
\[
T(m,h) \deq \pi_{m+1}(\inv{\M}\vect{p}^1(h_0))\,,
\]
where $\pi_{m+1}$ is the projection onto the $m+1$st coordinate.

Notice that 
\begin{equation}
    T(m,h) 
     = \begin{cases}
    \sum_{|z-w|=\ell^h} \av{h_0}{z}\av{h_0}{w} \neq 0 & \text{if } m = \ell^h \\
    0 & \text{if } m > \ell^h \,.
    \end{cases} \label{eq:chdist}
\end{equation}

By Proposition \ref{prop:sl}, there
exists for each $m$ and $j$ a sequence of random variables
$\{T_r(m,j)\}_{r=1}^\infty$ so that
almost surely,
\[
\lim_{r \rightarrow \infty} T_r(m,j) = T(m,h^j) \,.
\]
Let $\{s_n\}$ be any sequence with $\lim_{n \rightarrow \infty}s_n = \infty$.
An examination of the proof of Proposition \ref{prop:sl} shows that
there are bounds on the variances of $T_r(m,j)$ uniform for $m \leq
s_n$ and $j \leq s_n$. In particular, there is a sequence $f_r$ with
$\lim_{r \rightarrow \infty} f_r = 0$ so that
\[
\sup_{\substack{1 \leq j \leq s_n \\ \ell^j < m \leq s_n}}
\Eb{ T_r(m,j)^2} \leq f_r \,.
\]
 We have
\begin{align*}
\Pb{ \bigcup_{j=1}^{s_n} \bigcup_{m=\ell^j+1}^{s_n} \left\{ T_r(m,j) > \frac{1}{n} \right\}} 
& \leq s_n^2 n^2 \sup_{\substack{1 \leq j \leq s_n \\ \ell^j < m \leq s_n}}
\Eb{ T_r(m,j)^2} \\
& \leq s_n^2 n^2 f_r \,.
\end{align*}
Let $\{B_n\}$ be a sequence so that $s_n^2 n^2 f_{B_n}$ is summable.  By Borel-Cantelli, almost 
surely 
\[
T_{B_n}(m,j) \leq \frac{1}{n} \qquad \text{for } \ell^j < m \leq s_n \text{ and }
1 \leq j \leq s_n, \text{ eventually}\,.
\]
Similarly, since $\lim_{r \rightarrow \infty}T_r(\ell^j,j) =
T(\ell^j,j) \neq 0$, 
we can take a subsequence of $\{B_n\}$, which we will continue to denote
by $\{B_n\}$, so that almost surely
\[
T_{B_n}(\ell^j,j) > \frac{1}{n} \qquad \text{for } 1 \leq j \leq s_n,
\text{ eventually}\,.
\]
To summarize, the following holds almost surely: For $n$ sufficiently
large, for all $j=1,\ldots,s_n$,
\[
T_{B_n}(\ell^j,j) > n^{-1}, \text{ while }
T_{B_n}(m,j) \leq n^{-1} \text{ for } m =\ell^j +1,\ldots, s_n \,.
\]
Consequently,
\[
\sup\left\{ m \leq s_n \;:\; T_{B_n}(m,j) > \frac{1}{n} \right\} = \ell^j
\qquad \text{for } 1 \leq j \leq s_n \,,
\]
for $n$ sufficiently large.  Finally,
\[
\sup_{j \leq s_n} \sup\left\{ m \leq s_n 
\;:\; T_{B_n}(m,j) > \frac{1}{n} \right\} = \sup_{j \leq s_n} \ell^j
= \ell
\]
for $n$ sufficiently large. [Since $\sup_{j \leq r} \ell^j = \ell$ for
  all but finitely many $r$.]

\end{proof}  
We still need a proof of Proposition \ref{prop:sl}; it
is similar to the proof of \cite[Theorem 6.1]{LPP:PTRCT}.
Recall that $u_n = P^n(0,0)$.
The only facts we need about $P^n$ are
\begin{gather}
    \label{eq:f1}
    \sum_{n=0}^\infty u_n^2 = \infty \,,
    \intertext{and}
    \label{eq:f2}
    \lim_{n \rightarrow \infty} \frac{P^{n-k}(z,w)}{P^{n-j}(u,v)}
    = 1 \text{ for all } u,v,z,w \in \Z \text{ and } j,k \in \N \,.
\end{gather}
That the transition probabilities satisfy \eqref{eq:f1} and
\eqref{eq:f2} is
easily seen from a local central limit theorem, for example
\cite[Theorem 13.10]{W:RW}.
\proof[Proof of Proposition \ref{prop:sl}]
For 
\begin{equation*}
w(m,n)  \deq  \sum_{k=m+1}^n u_k^2  \quad
\text{and} \quad 
Z_k  \deq   \phi(X_{k})\phi(X_{k+t_1})\cdots \phi(X_{k+t_1 + \cdots + t_s}),
\end{equation*}
define the linear estimator
\begin{equation} \label{eq:lmndef}
L_{m,n}  \deq  \frac{1}{w(m,n)}
\sum_{k=m+1}^n u_k Z_k.
\end{equation}
Notice that
$\E[ L_{m,n} ]  =  \frac{1}{w(m,n)}
\sum_{k=m+1}^n u_k^2 
p_{\vect{t},k}(\phi)$,
where we define
\begin{equation*}
p_{\vect{t},k}  \deq  
\sum_{z\in \Z} \E[Z_k \given \rw_k = z ] \frac{P^k(0,z)}{u_k} \,.
\end{equation*}
Since
$\lim_{k \rightarrow \infty} \frac{u_k}{P^k(0,z)} = 1$
for all $z$,
it follows that $\lim_{k \rightarrow \infty}p_{\vect{t},k}(\phi) =
p_{\vect{t}}(\phi)$ and
\begin{equation} \label{eq:limitoflmn}
\lim_{m \rightarrow \infty}
\E[ L_{m,n} ]  =  p_{\vect{t}}(\phi).
\end{equation}

The strategy of the proof is to find sequences $\{m_i\}$ and $\{n_i\}$
so that
\begin{equation} \label{eq:whereto}
\lim_{N \rightarrow \infty} \frac{1}{N} \sum_{i=1}^N L_{m_i,n_i} 
 = p_{\bo{t}}\left( \phi \right) \text{ almost surely} \,.
\end{equation}
We prove (\ref{eq:whereto}) by showing that $\cov(Z_j,Z_k)$ is small
for $|j-k|$ large, and consequently $\cov(L_{m_i,n_i},L_{m_\tau,n_\tau})$ is small
for $|i-\tau|$ large.

For $\bo{z} \in \Z^{s+1}$ and $j \in \N$,
define the event
\begin{equation}
H_{\bo{t}}(j,\bo{z})   \deq   \{ \rw_{j} =
z_0 , \rw_{j+t_1} = z_1, \ldots, 
\rw_{j+t_1 + \cdots + t_s} = z_s\}.
\end{equation}
Then by conditioning on positions of $\rw$,
\begin{equation*}
\E[ Z_j Z_k ]  =  \sum_{\bo{z},\bo{w} \in \Z^{s+1}}
\E[ Z_j Z_k \given H_{\bo{t}}(j,\bo{z}) \cap 
H_{\bo{t}}(k,\bo{w}) ] \; \P[ H_{\bo{t}}(j,\bo{z}) \cap 
H_{\bo{t}}(k,\bo{w})]
\end{equation*}
Suppose that $k > j + t_1 + \cdots + t_s$.
Notice that
\begin{equation} \label{eq:ech}
\E[ Z_j Z_k \given H_{\bo{t}}(j,\bo{z}) \cap 
H_{\bo{t}}(k,\bo{w}) ] 
= 
\prod_{i=0}^s \av{\phi}{z_i}
\prod_{j=0}^s \av{\phi}{w_j},
\end{equation}
and
\begin{equation} \label{eq:ph}
\begin{split}
\P[ H_{\bo{t}}(j,\bo{z}) \cap 
H_{\bo{t}}(k,\bo{w})]
& =  P^j(0,z_0) \P[H_{\bo{t}}(j,\bo{z})  \given
S_j = z_0]  \\
& \quad \times P^{k-(j+t_1+\cdots+t_s)}(z_s,w_0)
 \P[H_{\bo{t}}(k,\bo{w}) \given S_k = w_0] \,.
\end{split}
\end{equation}
Combining (\ref{eq:ech}) and (\ref{eq:ph}) shows
that $\E[Z_j Z_k]$ is equal to
\begin{multline} \label{eq:ezjzk}
\sum_{\bo{z},\bo{w} \in \Z^{s+1}}
\left( \prod_{i=0}^s \av{\phi}{z_i}
\prod_{j=0}^s \av{\phi}{w_j} \right) 
P^j(0,z_0) \; \P[H_{\bo{t}}(j,\bo{z})  
\given S_j = z_0 ]  \\ 
\times
P^{k-(j+t_1+\cdots+t_s)}(z_s,w_0) \;
\P[H_{\bo{t}}(k,\bo{w}) \given S_k = w_0] \,.
\end{multline}

By similar reasoning, $\E[Z_j] \E[Z_k ] $ is equal to
\begin{multline*}
\sum_{\bo{z},\bo{w} \in \Z^{s+1}}
\left(
\prod_{i=0}^s \av{\phi}{z_i}
\prod_{j=0}^s \av{\phi}{w_j}
\right) \;
P^j(0,z_0) \; \P[H_{\bo{t}}(j,\bo{z})  
\given S_j = z_0]  \\ 
\times P^k(0,w_0) \;
\P[H_{\bo{t}}(k,\bo{w}) \given S_k = w_0] \,.
\end{multline*}

Because of the definitions of $a$ and $b$ in \eqref{eq:abdef},
if $D \geq |a| \vee |b|$, then for any $z \in \Z$ with
$|z| > D$, we have $\av{\phi}{z} = 0$.

Thus if $D \geq |a| \vee |b|$,
\begin{equation} \label{eq:covzjzk}
\begin{split}
|\cov(Z_j,Z_k)| & = \left| \; \E[Z_j Z_k ] - \E[Z_j] \E[Z_k ] \; \right|\\
& \leq
\sum_{\substack{-D \leq z_0 \leq D \\
-D \leq z_s \leq D\\ -D \leq w_0 \leq D}} P^j(0,z_0)
\left| P^{k-(j+t_1+\cdots+t_s)}(z_s,w_0) - P^k(0,w_0) \right|,
\end{split}
\end{equation}
since $\|\phi\|_{\infty} \leq 1$.

For sequences $\{m_i\}_{i=1}^\infty$ and 
$\{n_i\}_{i=1}^\infty$ with $m_i < n_i$, define
\begin{equation} \label{eq:lidef} 
L_i \deq L_{m_i,n_i} = \frac{1}{w(m_i,n_i)}
\sum_{j=m_i+1}^{n_i} u_j Z_j .
\end{equation}
For any sequence $\{\ep_i\}$, and integer $D$,
we will inductively define sequences $\{m_i\}$ and $\{n_i\}$ so that for
$\tau > i$,
\begin{equation} \label{eq:conclu}
 \left| \cov(L_i,L_\tau) \right| < C \ep_i \,, 
\end{equation}
provided $D \geq |a|\vee|b|$.  For now, assume
that $D$ is indeed chosen so that $D \geq |a| \vee |b|$.
To begin, let $m_1=n_1=1$.  Now assume that the pair $(m_i,n_i)$ has 
already been defined.
Pick $m_{i+1}$ large enough so that 
\begin{multline} \label{eq:pdiff}
    \text{if } j \leq n_i \text{ and } k>m_{i+1}, \text{then} \\
        \left| P^{k-(j+t_1+\cdots+t_s)}(z,w)-P^k(0,w) \right| \;
\leq \; C_1(D) \ep_i u_k 
        \text{ for all } -D \leq z,w \leq D.
\end{multline}
This is possible by \eqref{eq:f2}.

Next, pick $n_{i+1}$ so that
\begin{equation} \label{eq:nicond}
    w(m_{i+1},n_{i+1}) \geq w(m_{i+1})\,,
\end{equation}
where $w(m) \deq w(0,m)$.
Combining \eqref{eq:covzjzk} and \eqref{eq:pdiff} shows that
\begin{equation} \label{eq:covbnd}
        | \cov(Z_j,Z_k) |  \leq (2D)^3 \, C_1 (D) \, \ep_i \, u_j u_k 
        = C_2(D) \, \ep_i \, u_j u_k
\end{equation}
for $j<n_i$ and $k>m_{i+1}$.

{}From the definition \eqref{eq:lidef} and the bound \eqref{eq:covbnd} we have
\begin{equation*}
    \begin{split}
        \left| \cov(L_i,L_\tau) \right| 
& = \left| \frac{1}{w(m_i,n_i)}\frac{1}{w(m_\tau,n_\tau)}
        \sum_{j=m_i+1}^{n_i}
        \sum_{k=m_\tau+1}^{n_\tau} u_j u_k \cov(Z_j,Z_k) \right|\\
        & \leq \frac{1}{w(m_i,n_i)}\frac{1}{w(m_\tau,n_\tau)}
        C_2(D) \ep_i \sum_{j=m_i+1}^{n_i}u_j^2 \sum_{k=m_\tau+1}^{n_\tau} 
        u_k^2 \\
        & = C_2(D) \ep_i \,.
    \end{split}
\end{equation*}

Next, we show that
\begin{equation} \label{eq:boundedsecondmoment}
\sup_{i \geq 1} \E[L_i^2] < M < \infty.
\end{equation}
Notice that
\begin{equation} \label{eq:lisquared}
\E[ L_i^2 ]  \leq  \frac{1}{w(m_i,n_i)^2}
\left\{ 2\sum_{\substack{m_i < 
j\leq k \leq n_i \\ |j-k|\leq t_1 + \cdots + t_s}}
u_ju_k\E[Z_jZ_k] +
2\sum_{\substack{m_i < j,k \leq n_i \\ j+t_1+\cdots+t_s < k }}u_ju_k \E[Z_jZ_k]
\right\}
\end{equation}
By \eqref{eq:f2},  for each fixed $k$, we have
$\lim_{j \rightarrow \infty}\frac{u_j}{u_{j+k}} = 1$,
and hence
$u_{j+k} \leq C_3 u_j$ for $k=0,1,\ldots,t_1+\cdots+t_s$.  Consequently,
since $\| \phi \|_{\infty} \leq 1$ and hence $Z_jZ_k \leq 1$,
\begin{equation} \label{eq:close}
\sum_{\substack{m_i < j\leq k \leq n_i \\ |j-k|\leq t_1+\cdots+t_s}}
u_ju_k\E[Z_jZ_k] 
 \leq  
(t_1+\cdots+t_s) C_{3} \sum_{j=m_i+1}^{n_i} u_j^2  = 
C_{4} w(m_i,n_i) \,.
\end{equation}
Since $\E[Z_j Z_k]$ is equal to \eqref{eq:ezjzk}, we have
\begin{equation*} 
\E[Z_j Z_k]
 \leq  \sum_{z_0=-D}^{D}
P^j(0,z_0) \sum_{z_s=-D}^{D} \sum_{w_0=-D}^{D} P^{k-(j+t_1+\cdots+t_s)}(z_s,w_0)
 \leq  C_{5}(D) u_j u_{k-j},
\end{equation*}
again by \eqref{eq:f2}. 
Thus,
\begin{equation} \label{eq:farapart}
2\sum_{\substack{m_i < j,k \leq n_i \\ j+t_1+\cdots+t_s < k}}u_ju_k \E[Z_jZ_k]
 \leq 
C_{6}(D) \sum_{j=m_i+1}^{n_i} u_j^2 \sum_{k=j+1}^{n_i}
u_k u_{k-j} \,.
\end{equation}
Now applying the Cauchy-Schwarz inequality yields
\begin{equation} \label{eq:anothercs}
\sum_{k=j+1}^{n_i} u_k u_{k-j}
\leq  \sqrt{w(j,n_i)w(1,n_i - j)}  \leq  w(n_i) .
\end{equation}
Plugging (\ref{eq:anothercs}) into (\ref{eq:farapart}) yields
\begin{equation} \label{eq:farapart2}
2\sum_{\substack{m_i < j,k \leq n_i \\ j+t_1+\cdots + t_s<k}}u_ju_k \E[Z_jZ_k]
 \leq  C_{6}(D) w(n_i) w(m_i,n_i) \,.
\end{equation}
Then, using (\ref{eq:farapart2}) and (\ref{eq:close}) in 
(\ref{eq:lisquared}) yields
\begin{equation*}
\E[ L_i^2 ]  \leq  
\frac{C_{4}}{w(m_i,n_i)} + C_{6}(D) \frac{w(n_i)}{w(m_i,n_i)} \,.
\end{equation*}
Now, since by \eqref{eq:nicond} 
$w(m_i,n_i) \geq w(m_i)$, it follows that $w(n_i) \leq 2 
w(m_i,n_i)$.  Hence,
\begin{equation*}
\E[ L_i^2 ]  \leq  
\frac{C_4}{w(m_i,n_i)} + {2C_6(D)}  \leq  M \,,
\end{equation*}
where $M$ doesn't depend on $i$.

If we choose $\ep_i = i^{-3}$, then a strong law for weakly 
correlated and centered random variables (see \cite[Theorem 37.7.A]{L:PT2})
applied to the sequence
$\{L_i - \E[L_i]\}_{i \geq 1}$ implies that
$\frac{1}{N}\sum_{i=1}^N \left( L_i - \E[L_i] \right)
\rightarrow 0$.

Thus, from (\ref{eq:limitoflmn}), we have 
$\lim_{i \rightarrow \infty}\E[L_i]  =
 p_{\vect{t}}(\phi)$
and 
$\frac{1}{N}\sum_{i=1}^N L_i \rightarrow 
p_{\vect{t}}(\phi)$.

To summarize, for each $D$ we get a sequence of
functions $\left\{g_N^D\right\}_{N=1}^\infty$ so that
\begin{equation} \label{eq:glimit}
  \text{if } D \geq |a|\vee |b|, \text{ then }
    g_N^{D}(X_1,\ldots,X_N) \xrightarrow{N \rightarrow \infty} 
    p_{\vect{t}}(\phi) \text{ a.s.} 
\end{equation}
By a diagonalization argument, there is a sequence of
integers $\{r_n\}$
so that for each $D$, the limit $A(D) \deq \lim_{n} g_{r_n}^D
(X_1,\ldots,X_{r_n})$ exists, and moreover
for each $D$, 
\[
|g_{r_n}^D(X_1,\ldots,X_{r_n}) - A(D)| < n^{-1} \text{ for } 
n \geq D \,.
\]
From \eqref{eq:glimit} it follows that
$A(D) = p_{\vect{t}}(\phi)$ for $D \geq |a|\vee |b|$, and consequently a.s.
\begin{equation*}
    g^{D}_{r_{D}}(X_1,\ldots,X_{r_{D}}) \xrightarrow{D \rightarrow \infty} 
    p_{\vect{t}}(\phi) \,.
\end{equation*}
\endproof

\section{Proofs of Theorem \ref{thm:coins} and
Theorem \ref{thm:main}} \label{sec:proof}
Here is an outline of the proof of Theorem \ref{thm:main}:
In the last section, it was shown in Proposition \ref{prop:sl} that 
the parameters $\{p_{\bo{t}}\}$ can be recovered 
{}from the data $\{X_1,X_2,\ldots \}$.  We transformed these
parameters in Proposition \ref{prop:ptoq} into the
family $\{Q_{\bo{d}}\}$.  Now we must show how to use
this information to determine the laws $\{\eta(z)\}$.

Before we prove Theorem \ref{thm:main} in its full generality,
we offer the proof of Theorem \ref{thm:coins},
the case where the observable random variables are coin tosses.
The exposition is much clearer in
this special case.  The proof is given in its full generality 
later in this section.

\medskip
\noindent
{\bf Theorem \ref{thm:coins}.} \em
Let $S$ be a symmetric and aperiodic random walk, and 
label each $z \in \Z$ with a coin having bias $\bias(z)$.
If only finitely many vertices have $\bias(z)\neq 0$,
then there exist a sequence $\{\bias_n\}_{n=0}^\infty$ with
$\bias_n$ measurable with respect to $\sigma(X_1,\ldots,X_n)$,
and a
shift and reflection $T$ (depending on $\eta$) so that
$$
\lim_{n \rightarrow \infty} \theta_n(z) = \theta \circ T(z) \text{ for all }
z \in \Z\,, \text{ a.s.}
$$
\rm

Let us first introduce some notation.
\begin{defn} \label{def:symdef}
    Let $(x_1,x_2,\ldots,x_n)$ be a finite, ordered sequence.
    Define $[x_1,x_2,\ldots,x_n]$ as
    \begin{equation} \label{eq:symdef}
        [x_1,x_2,\ldots,x_n] \deq
        \left\{ (x_1,x_2,\ldots,x_n), (x_n,x_{n-1},\ldots,x_1) \right\} \,.
    \end{equation}
\end{defn}

\begin{prop} \label{prop:needed}
    Let $\theta:\Z \rightarrow [0,1]$ have compact support.
    Then, given $\ell=b-a$ as defined in \eqref{eq:abdef},
    \begin{equation} \label{eq:bias_equiv}
    [\bias(a),\bias(a+1),\ldots,\bias(b-1),\bias(b)] 
    \end{equation}
    can be uniquely determined from $\{\vect{p}^k(I)\}$,
    where $I$ is the identity map.
\end{prop}

\begin{proof}
We will write $Q^k_{\vect{d}}$ for $Q^k_{\vect{d}}(I)$,
where $I$ is the identity map.

By Proposition \ref{prop:ptoq} 
it is enough to construct
\eqref{eq:bias_equiv} using only $\{\vect{Q}^k\}$.
    
Our method will be to first determine the biases
at $a$ and $b$, and then work inwards.  

In the case where $k=1$ and $\vect{d} = (\ell)$, we have
\begin{equation*}
\frac{1}{2} {Q^1_\ell} = \frac{1}{2}\sum_{\substack{(z_1,z_2) \\ 
|z_2-z_1| = {\ell}}}
\bias(z_1)\bias(z_2) = \bias(a)\bias(b) \,,
\end{equation*}
as
the only pairs $(z_1,z_2)$ with 
$|z_2-z_1| = {\ell}$ and with
$\bias(z_1)\bias(z_2) \neq 0$ are $(a,b)$ and $(b,a)$.

Notice that
\begin{equation*}
    Q^2_{(\ell,\ell)} =
    \sum_{\substack{(z_1,z_2,z_3) \\ |z_2-z_1| = \ell, \\
    |z_3 - z_2| = \ell}} \bias(z_1)\bias(z_2)\bias(z_3) 
    = \bias(a)\bias(b)\left\{ \bias(a) + \bias(b) \right\} \,,
\end{equation*}
and so
$\bias(a) + \bias(b) =
2 Q^2_{(\ell,\ell)}/Q^1_{\ell}$. 
We can conclude that
\begin{align*}
    \{ \bias(a) , \bias(b) \} & =
    \left\{ \frac{(\bias(a)+\bias(b)) \pm 
    \sqrt{(\bias(a)+\bias(b))^2 - 4\bias(a)\bias(b)}}{2} \right\} \\
    &= \left\{   Q^2_{(\ell,\ell)}/ Q^1_{\ell} \pm
    \sqrt{ \left( Q^2_{(\ell,\ell)}/ Q^1_{\ell}\right)^{2} 
    -Q^{1}_{\ell}/2 } \right\}\,.
\end{align*}
    
Now we move on to the unordered pair $\{\bias(a+1),\bias(b-1)\}$.
Notice that
\begin{equation*}
    Q^2_{(\ell,1)}  = \sum_{\substack{(z_1,z_2,z_3) \\ |z_2-z_1| = 
    \ell, \\
    |z_3 - z_2| = 1}} \bias(z_1)\bias(z_2)\bias(z_3) 
    = \bias(a)\bias(b)\left\{ \bias(a+1) + \bias(b-1) \right\} \,,
\end{equation*}
and
\begin{align*}
    Q^5_{(\ell,1,1,\ell,1)} & =\sum_{\substack{(z_1,z_2,z_3,z_4,z_5)
    \\ |z_2-z_1| = \ell, 
    |z_3 - z_2| = 1, \\ |z_4 - z_3| = 1, |z_5 - z_4| = \ell \\
    |z_6-z_5|=1}} 
    \bias(z_1)\bias(z_2)\bias(z_3)\bias(z_4)\bias(z_5) \bias(z_6)\\
    & = 2\left\{ \bias(a)\bias(b) \right\}^2 \bias(a+1) \bias(b-1)
    \,.
\end{align*}
Hence,
\begin{equation*}
\bias(a+1) + \bias(b-1) = 2Q^{2}_{(\ell,1)}/Q^{1}_{\ell} \quad
\text{ and } \quad
\bias(a+1)\bias(b-1) = 2Q^{5}_{(\ell,1,1,\ell,1)}/
\left( Q^{1}_{\ell} \right)^{2} \,,
\end{equation*}
and so the unordered pair $\{\bias(a+1),\bias(b-1)\}$ can be
written as a function of $(\vect{Q}^1,\ldots,\vect{Q}^5)$.

The cross-product
$\bias(a)\bias(a+1) + \bias(b)\bias(b+1)$
allows us to write $[\bias(a),\bias(a+1),
\bias(b-1),\bias(b)]$ as a function of the unordered pairs
$\{\bias(a),\bias(b)\}$ and $\{\bias(a+1),\bias(b-1)\}$.
We have
\begin{align*}
    Q^3_{(\ell,\ell,1)} & = \sum_{\substack{(z_1,z_2,z_3,z_4) \\
    |z_2 - z_1| = \ell, |z_3 - z_2| = \ell, \\
    |z_4 - z_3| = 1}} \bias(z_1)\bias(z_2)\bias(z_3)\bias(z_4) \\
    & = \bias(a)\bias(b) \; \left\{ \bias(a)\bias(a+1) + 
    \bias(b)\bias(b11) \right\} \,,
\end{align*}
and so 
$$
\bias(a)\bias(a+1) + \bias(b)\bias(b-1) =
2Q^3_{(\ell,\ell,1)}/Q^1_\ell \,.
$$
We conclude that $[\bias(a),\bias(a+1),\bias(b-1),\bias(b)]$
can be written as a function of the $\{\vect{Q}^{k}\}$.

Continuing in this way, we work down to the center to obtain
\eqref{eq:bias_equiv},
using only $\{\vect{Q}^{k}\}$.
\end{proof}

\begin{proof}[Proof of Theorem \ref{thm:coins}]
    If $\ell$ is given, then by Proposition \ref{prop:needed}, it is
    clear that the theorem holds.  That is, for each $\ell$, we can
    construct a scenery $\theta^\ell_n$ depending only on
    $\sigma(X_1,\ldots,X_n)$ measurable random variables,
    and so that $\theta^\ell_n \rightarrow \theta\circ T$.
    Then by Lemma \ref{lem:balem}, the sequence
    $\theta^{\ell_n}_n$ will satisfy the requirements of the 
    Theorem.
\end{proof}

We now provide a proof of our result in its full generality.
Recall the notation given in Definition \ref{defn:func_ave}:
As $\eta(z)$ is a probability measure on $\R$ for each $z\in\Z$,
$$
\av{\phi}{z} = \int_\R \phi d\eta(z)\,, 
\quad \text{and} \quad
\phi_0 = \phi - \int_\R \phi d\alpha \,.
$$

\proof[Proof of Theorem \ref{thm:main}]
We will show how to use the collection $\Qs \deq \{\vect{Q}^k(\phi)\}_{\phi 
\in \Phi_1,k\in \N}$ and $\ell$ to construct a stochastic scenery
which is a shift and/or a reflection of $\eta$.  The conclusion
of the Theorem will then follow from applications of
Proposition \ref{prop:sl}, Proposition \ref{prop:ptoq}, and
Lemma \ref{lem:balem} similarly to the proof of
Theorem \ref{thm:coins}.

Fix a bounded and measurable $\phi$.  
We have that
\begin{equation} \label{eq:st1}
\begin{split}
 Q^1_\ell(\phi_0)
& = 
2 \av{\phi_0}{a}\av{\phi_0}{b} \,, \\
Q^2_{(\ell,\ell)}(\phi_0) & = 
\av{\phi_0}{a}\av{\phi_0}{b}\left\{ \av{\phi_0}{a} + 
\av{\phi_0}{b}\right\} \,.
\end{split}
\end{equation}
Assume first that $\av{\phi_0}{a}\av{\phi_0}{b} \neq 0$.  
In this case, we can solve in \eqref{eq:st1} to obtain
\begin{equation} \label{eq:st2}
        \av{\phi_0}{a}\av{\phi_0}{b} =
        Q^1_\ell(\phi_0) /2 \quad \text{ and } \quad
        \av{\phi_0}{a} + \av{\phi_0}{b} =
        2 Q^2_{(\ell,\ell)}(\phi_0)/Q^1_\ell(\phi_0) \,.
\end{equation}
We use the identities in \eqref{eq:st2} to get
\begin{equation} \label{eq:st2a}
\{\av{\phi_0}{a}, \av{\phi_0}{b}\}
=
\left\{ \frac{Q^2_{(\ell,\ell)}(\phi_0)}
{Q^1_\ell(\phi_0)} \pm
\sqrt{ \left(\frac{Q^2_{(\ell,\ell)}(\phi_0)}{Q^1_\ell(\phi_0)}
\right)^2 - Q^1_\ell(\phi_0)/2 }
\right\}
\end{equation}

We proceed inductively: 
assume we have determined, for $k \leq \ell/2$,
\begin{equation} \label{eq:indhypr}
\{ (\av{\phi_0}{a},\ldots,\av{\phi_0}{a+k-1}), \,
(\av{\phi_0}{b},\ldots,\av{\phi_0}{b-k+1}) \} \,.
\end{equation}
We have
\begin{equation} \label{eq:st3}
    \begin{split}
        Q^4_{(\ell,k,k,\ell)}(\phi_0) 
        & = \left(\av{\phi_0}{a}\av{\phi_0}{b} \right)^2 
        \left\{ \av{\phi_0}{a+k} + \av{\phi_0}{b-k} \right\} \,, \\
        Q^5_{(\ell,k,k,\ell,k)}(\phi_0)
        & = 2\left(\av{\phi_0}{a}\av{\phi_0}{b} \right)^2
        \av{\phi_0}{a+k}\av{\phi_0}{b-k}\,.
    \end{split}
\end{equation}
By assumption,
$\av{\phi_0}{a}\av{\phi_0}{b} \neq 0$, so 
we can solve in \eqref{eq:st3} (using \eqref{eq:st2}) to get
\begin{equation*}
        \av{\phi_0}{a+k}\av{\phi_0}{b-k}
         = \frac{2Q^5_{(\ell,k,k,\ell,k)}}
        {\left( Q^1_\ell(\phi_0) \right)^2} \quad \text{and} \quad
        \av{\phi_0}{a+k} + \av{\phi_0}{b-k}
         = \frac{4Q^4_{(\ell,k,k,\ell)}(\phi_0) }
        {\left( Q^1_\ell(\phi_0) \right)^2} \,.
\end{equation*}
Thus, similarly to \eqref{eq:st2a}, we can determine
\begin{equation} \label{eq:onemore}
\{ \av{\phi_0}{a+k},\av{\phi_0}{b-k} \} 
\end{equation}
as a function of elements of $\Qs$.

Now define
$m \deq \sup \{j \leq k-1 : \av{\phi_0}{a+j} \neq \av{\phi_0}{b-j} 
\}$.
If there is no such $j$, there is no problem in finding
\begin{equation} \label{eq:indconcr}
\{ (\av{\phi_0}{a},\ldots,\av{\phi_0}{a+k}), \,
(\av{\phi_0}{b},\ldots,\av{\phi_0}{b-k}) \} 
\end{equation}
{}from \eqref{eq:indhypr} and \eqref{eq:onemore},
so without loss of generality, assume such a $j$ exists.

Now notice that
\begin{gather} 
    Q^7_{(\ell,m,m,\ell,k,k,\ell)} = 
    \left(\av{\phi_0}{a}\av{\phi_0}{b} \right)^3
    \left\{ \av{\phi_0}{a+m}\av{\phi_0}{b-k} + 
    \av{\phi_0}{b-m}\av{\phi_0}{a+k} \right\} \nonumber \,, \\
     \intertext{and hence}
    \av{\phi_0}{a+m}\av{\phi_0}{b-k} + 
    \av{\phi_0}{b-m}\av{\phi_0}{a+k} 
    = \frac{8Q^7_{(\ell,m,m,\ell,k,k,\ell)} }
    { \left( Q^1_\ell(\phi_0) \right)^3   } \,.
    \label{eq:cp}
\end{gather}
{}From the two unordered pairs (which we already know)
\begin{gather*}
    \{ \av{\phi_0}{a+m},\av{\phi_0}{b-m} \} \quad
    \text{ and } \quad
    \{ \av{\phi_0}{a+k},\av{\phi_0}{b-k} \}
\end{gather*}
we can get two possible pairings:
\begin{equation} \label{eq:twopair}
    \begin{array}{c}
F_{good} \deq \{ \, (\av{\phi_0}{a+m},\av{\phi_0}{a+k}), \, 
(\av{\phi_0}{b-m},\av{\phi_0}{b-k}) \, \} \\
\text{and} \\
F_{bad} \deq \{ \, (\av{\phi_0}{a+m},\av{\phi_0}{b-k}), \,
(\av{\phi_0}{b-m},\av{\phi_0}{a+k}) \, \} \,.
\end{array}
\end{equation}

{}From these two pairings we can compute the values
\begin{equation} \label{eq:twovals}
\begin{array}{c}
f_{good} \deq \av{\phi_0}{a+m}\av{\phi_0}{a+k} +
\av{\phi_0}{b-m}\av{\phi_0}{b-k} \\
\mbox{ and } \\
f_{bad} \deq \av{\phi_0}{a+m}\av{\phi_0}{b-k} +
\av{\phi_0}{b-m}\av{\phi_0}{a+k}
\end{array} \,.
\end{equation}
Since we have written $f_{bad}$ in 
\eqref{eq:cp} as a function of elements of $\Qs$,
we can use these elements to discard the pairing $F_{bad}$
in (\ref{eq:twopair}) 
{}from which it comes.  If $f_{bad} =
f_{good}$, we cannot determine the 
proper pairing. But this happens if and only if
\begin{equation} \label{eq:whym}
\av{\phi_0}{a+k}\left( \av{\phi_0}{a+m} - \av{\phi_0}{b-m} \right)
= \av{\phi_0}{b-k}\left( \av{\phi_0}{a+m} - \av{\phi_0}{b-m} 
\right) .
\end{equation}
Since $m$ was defined so that 
$\av{\phi_0}{a+m} - \av{\phi_0}{b-m} \neq 0$,
we
conclude that (\ref{eq:whym}) is valid if and only if
$\av{\phi_0}{a+k} = \av{\phi_0}{b-k}$.  But in this case, it is 
clear that (\ref{eq:indconcr}) can be found {}from
\eqref{eq:indhypr} and \eqref{eq:onemore}.

Otherwise, we can pick out $F_{good}$ from
(\ref{eq:twopair}), and hence we can determine (\ref{eq:indconcr}).

Thus, by induction, we can determine
\begin{equation} \label{eq:recfcn}
\{ (\av{\phi}{a},\ldots,\av{\phi}{b}),
(\av{\phi}{b},\ldots,\av{\phi}{a}) \} \,.
\end{equation}

To summarize, we have just shown for any bounded measurable function 
$\phi$ obeying the condition
$\av{\phi_0}{a}\av{\phi_0}{b} \neq 0$ we can determine
\begin{equation} \label{eq:uptorev}
\lbr \av{\phi}{a},\ldots,\av{\phi}{b} \rbr \,.
\end{equation}

Now suppose that $\av{\phi_0}{a}\av{\phi_0}{b} = 0$.
Recall that the elements of $\Phi_1$ are enumerated
as $(h_1,\ldots,h_n,\ldots)$.  If we define
$n^*$ as the smallest $n$ so that 
$\{ \av{h_n}{a}-\int_\R h_n d\al ,\av{h_n}{b} -\int_\R h_n d\al \}$
has distinct non-zero elements, then Lemma \ref{lem:disting}
guarantees that $n^* < \infty$.
Let $\psi \deq h_{n^*}$.

There 
is a $c \in \Q$ small enough so that
\begin{equation*}
\av{\phi_0 + c\psi_0}{a} \neq 0 \quad \mbox{and} \quad
\av{\phi_0 + c \psi_0}{b} \neq 0 \,.
\end{equation*}
We can determine
    $\lbr \av{c\psi_0}{a},\ldots, \av{c\psi_0}{b} \rbr$,
and 
    $\lbr \av{\phi_0 + c \psi_0}{a},\ldots,\av{\phi_0 + c \psi_0}{b} \rbr$.

Thus we can determine the two quantities
\begin{equation*}
    \begin{array}{c}
\lbr \av{\phi_0 + c \psi_0}{a} - \av{c\psi_0}{a},\ldots,\av{\phi_0 + c 
\psi_0}{b} - 
\av{c\psi_0}{b} \rbr\\
\mbox{and}\\
\lbr \av{\phi_0 + c \psi_0}{a} - \av{c\psi_0}{b},\ldots,\av{\phi_0 + c 
\psi_0}{b} - 
\av{c\psi_0}{a} \rbr\\
\end{array}.
\end{equation*}
But the end points of the second expression are nonzero for $c$ small 
enough, while at least one endpoint in the first express is zero.
Thus we can pick out the first expression, which allows us to 
determine (\ref{eq:uptorev}).

The careful reader will have noticed that our use of
Lemma \ref{lem:disting}
requires that $\eta(a) \neq \eta(b)$.  If this is not the
case, then $\av{\phi_0}{a} = \av{\phi_0}{b} = 0$.  We can
then apply the entire procedure with $\ell-2,a+1,b-1$ replacing
$\ell,a$ and $b$ respectively.

Now we have shown that we can determine
\begin{equation} \label{eq:firstinfo}
\left(\, \lbr \av{\phi_1}{a},\ldots,\av{\phi_1}{b} \rbr, \,
\lbr \av{\phi_2}{a},\ldots,\av{\phi_2}{b} \rbr, \ldots \, \right) \,.
\end{equation}
Our goal now is to find
\begin{equation*}
    \lbr \left( \av{\phi_k}{a} \right)_{k=1}^\infty,
    \left( \av{\phi_k}{a+1} \right)_{k=1}^\infty,
    \ldots, \left( \av{\phi_k}{b} \right)_{k=1}^\infty
    \rbr \,,
\end{equation*}
which, since $\Phi$ is measure determining, yields
    $\lbr \eta(a), \eta(a+1), \ldots, \eta(b) \rbr$.
Again, we use induction. 
Suppose that we know
\begin{equation} \label{eq:indhypgen}
\lbr \left(\, \av{\phi_1}{a},\ldots,\av{\phi_{k-1}}{a} \right), \,
\ldots,\, \left( \av{\phi_1}{b},\ldots,\av{\phi_{k-1}}{b} \right) \, \rbr \,.
\end{equation}
{}From (\ref{eq:firstinfo}), we have in particular
\begin{equation} \label{eq:kinfo}
\lbr \av{\phi_k}{a},\ldots,\av{\phi_k}{b} \rbr \,.
\end{equation}

There are two possible extensions of (\ref{eq:indhypgen}) using 
(\ref{eq:kinfo}):
\begin{equation} \label{eq:twoextens}
\begin{array}{c}
    G_{good} \deq \lbr \left(\av{\phi_1}{a},\ldots,\av{\phi_{k-1}}{a},\av{\phi_k}{a}
    \right), \ldots,  
    \left(\av{\phi_1}{b},\ldots,\av{\phi_{k-1}}{b},\av{\phi_k}{b} \right)
    \rbr \\
    \mbox{and}\\
    G_{bad} \deq \lbr \left(\av{\phi_1}{a},\ldots,\av{\phi_{k-1}}{a},\av{\phi_k}{b}
    \right), \ldots,  
    \left(\av{\phi_1}{b},\ldots,\av{\phi_{k-1}}{b},\av{\phi_k}{a} \right)
    \rbr 
\end{array}
\end{equation}
We want to guarantee that we can pick out $G_{good}$ from 
$\{G_{good},G_{bad} \}$ given in 
(\ref{eq:twoextens}).

Suppose there is some $j$ so that
\begin{equation} \label{eq:difj}
\left( \av{\phi_1}{a+j},\ldots,\av{\phi_{k-1}}{a+j}\right) \; \neq \;
\left( \av{\phi_1}{b-j},\ldots,\av{\phi_{k-1}}{b-j}\right) \,.
\end{equation}
This is without loss of generality, since
if there is no such $j$, then of course we can extend 
(\ref{eq:indhypgen}) to obtain
\begin{equation} \label{eq:indgen}
\lbr \left( \av{\phi_1}{a},\ldots,\av{\phi_{k}}{a} \right),
\ldots, \left( \av{\phi_1}{b},\ldots,\av{\phi_{k}}{b} \right) \rbr \,.
\end{equation}
Since the linear functionals on 
$\R^{k-1}$ separate points, there exist $\bo{c} \in \R^{k-1}$ so that
if $T_c(\bo{x}) = \langle \bo{c}, \bo{x} \rangle$, then
\begin{equation*}
T_c\left( \, (\av{\phi_1}{a+j},\ldots,\av{\phi_{k-1}}{a+j}) \, \right)
\neq T_c\left(\, (\av{\phi_1}{b-j},\ldots,\av{\phi_{k-1}}{b-j}) \, \right)
\,.
\end{equation*}
Since for any function $\phi$ we can recover (\ref{eq:recfcn}), by 
linearity of expectation, for the function 
$T_c \left( \, (\phi_1,\ldots,\phi_{k-1}) \, \right)
+ \phi_k$
we can determine
\begin{equation} \label{eq:det}
\lbr T_c( \av{\phi_1}{a},\ldots,\av{\phi_{k-1}}{a}) +\av{\phi_k}{a},\ldots,
T_c( \av{\phi_1}{b},\ldots,\av{\phi_{k-1}}{b}) +\av{\phi_k}{b} ) \rbr
\end{equation}
But from (\ref{eq:twoextens}), we can determine
\begin{equation}\label{eq:tdet}
\begin{array}{c}
    \lbr T_c\left(\av{\phi_1}{a},\ldots,\av{\phi_{k-1}}{a}\right)
    + \av{\phi_k}{a}, \ldots,  
    T_c\left(\av{\phi_1}{b},\ldots,\av{\phi_{k-1}}{b}\right)
    + \av{\phi_k}{b} \rbr \\
    \mbox{and}\\
    \lbr T_c\left(\av{\phi_1}{a},\ldots,\av{\phi_{k-1}}{a}\right)
    + \av{\phi_k}{b}, \ldots,  
    T_c\left(\av{\phi_1}{b},\ldots,\av{\phi_{k-1}}{b}\right) + \av{\phi_k}{a} 
    \rbr 
\end{array}
\end{equation}
It thus suffices to show that these two sets cannot be the same.  
If the sets in (\ref{eq:tdet}) are the same, then
either $\av{\phi_k}{a+j} = \av{\phi_k}{b-j}$ for all $j$, in which 
case there is no problem in obtaining (\ref{eq:indgen}), or
\begin{equation*}
T_c\left( \av{\phi_1}{a+j},\ldots,\av{\phi_{k-1}}{a+j}\right) \; = \;
T_c\left( \av{\phi_1}{b-j},\ldots,\av{\phi_{k-1}}{b-j}\right) \;
\mbox{ for all } j \,.
\end{equation*}
But $T_c$ was defined to rule out this possibility.  Hence we can
distinguish the ``correct'' extension in (\ref{eq:twoextens}), and thus
we can determine (\ref{eq:indgen}).

By induction, we can then determine
\begin{equation} \label{eq:ldver}
[ (\av{\phi_k}{a})_{k=1}^\infty,\ldots,(\av{\phi_k}{b})_{k=1}^\infty] 
\,.
\end{equation}
Because the $\{ \phi_k\}$ are measure determining, from
(\ref{eq:ldver}) we can determine
\begin{equation*}
[ \eta(a),\ldots,\eta(b)] \,.
\end{equation*}
Pick one $\ell$-tuple from $[ \eta(a),\ldots,\eta(b)]$, and
call it $(\eta_0(1),\ldots,\eta_0(\ell))$.
Finally, define $\hat{\eta}$ by
\begin{equation*}
\hat{\eta}(z) \; = \;
\left\{
\begin{array}{ll}
    \alpha & \mbox{if } z < 1 \mbox{ or } z > \ell \\
    \eta_0(z) & \mbox{if } 1 \leq z \leq \ell
\end{array}
\right. \,.
\end{equation*}
\endproof

\section{Unsolved Problems} \label{sec:conclusion}
There remain a couple of questions still requiring resolution:
\begin{enumerate}
\item
It follows from results in \cite{L:IS} that 
if we allow infinitely many biased coins, then there 
exist bias configurations  $\{\theta(z) : z \in \Z\}$ that cannot be
recovered up to a reflection and shift. We have restricted attention
to finitely many biased coins; however, the methods
of this paper do allow recovery of infinitely many biases
(up to reflection and shift) if the biased coins are sufficiently sparse.
It remains an  open question to determine for which
infinite sets of integers $S$  can one recover the biases
of biased coins located in $S$ from the tosses observed by a random walker.  
\item
We have shown in Theorem \ref{thm:main} 
that for {\em symmetric} walks with finite variance increments, we can
recover finitely many unknown distributions up to a shift and/or a reflection.
We have not been able to generalize the algebraic results in
Section \ref{sec:propositions} to walks which are asymmetric,
yet still mean zero.  
In particular, Lemma \ref{lem:fullrank} needs to be
proven for a matrix with entries $\{P^i(0,j) \,:\,
1 \leq i \leq r, \ -m \leq j \leq m\}$.
We expect the following to hold:
\begin{conj}
In the setting of theorem \ref{thm:main}, suppose that the symmetry assumption
on the walk $S$ is replaced by the assumption
that the walk $S$ has mean-zero, non-symmetric increments of finite
variance. Then $\eta$ can be recovered up to a shift only.
\end{conj}

\item In Theorem \ref{thm:main}, is it possible to replace 
the symmetric random walk on $\Z$ satisfying (\ref{eq:qconds1}), 
by a general null-recurrent 
Markov chain that satisfies $\sum_n u_n^2=\infty$? \newline
(Here $u_n$ is the probability that the chain returns to  its starting state
 after $n$ steps.) The goal is to recover $\eta$ up to an automorphism
of the Markov chain.
\end{enumerate}

\medskip
\noindent
{\em Acknowledgements.}
We are grateful  to Rich Bass for useful discussions,
and to Alan Hammond and Gabor Pete for helpful comments on 
the manuscript.


\begin{thebibliography}{How96b}

\bibitem[BK96]{BK:Scen}
I.~Benjamini and H.~Kesten.
\newblock Distinguishing sceneries by observing the scenery along a random walk
  path.
\newblock {\em J.\ Anal.\ Math.}, 69:97--135, 1996.

\bibitem[HJ91]{HJ:TMA}
R.A. Horn and C.R. Johnson.
\newblock {\em Topics in matrix analysis}.
\newblock Cambridge University Press, Cambridge, 1991.

\bibitem[HK97]{HK:RCT}
M.~Harris and M.~Keane.
\newblock Random coin tossing.
\newblock {\em Probab.\ Th.\ Rel.\ Fields}, 109:27--37, 1997.

\bibitem[How96a]{H:PS}
C.D. Howard.
\newblock Detecting defects in periodic scenery by random walks on {${\bf Z}$}.
\newblock {\em Rand.\ Struc.\ Algorithms}, 8(1):59--74, 1996.

\bibitem[How96b]{H:RWS}
C.D. Howard.
\newblock Orthogonality of measures induced by random walks with scenery.
\newblock {\em Combin.\ Probab.\ Comput.}, 5:247--256, 1996.

\bibitem[Kes98]{K:Scenery}
H.~Kesten.
\newblock Distinguishing and reconstructing sceneries from observations along
  random walk paths.
\newblock In D.~Aldous and J.~Propp, editors, {\em Microsurveys in Discrete
  Probability}, volume~41 of {\em DIMACS Series in Discrete Mathematics and
  Theoretical Computer Science}. American Mathematical Society, 1998.

\bibitem[Lin99]{L:IS}
E.~Lindenstrauss.
\newblock Indistinguishable sceneries.
\newblock {\em Rand.\ Struc.\ Algorithms}, 14(1):71--86, 1999.

\bibitem[Lo{\`e}78]{L:PT2}
M.~Lo{\`e}ve.
\newblock {\em Probability Theory II}.
\newblock Springer-Verlag, New York, NY, fourth edition, 1978.

\bibitem[LM02]{LM:2d}
M.~L\"owe and H.~Matzinger.
\newblock Scenery reconstruction in two
dimensions with many colors. 
\newblock {\em Ann.\ Appl.\ Probab.}
12(4): 1322--1347, 2002.

\bibitem[LM03]{LM:Corr}
M.~L\"owe and H.~Matzinger.
\newblock Reconstruction of sceneries with
correlated colors. 
\newblock {\em Stochastic Process.\ Appl.} 
105(2): 175--210, 2003. 

\bibitem[LPP01]{LPP:PTRCT}
D.~Levin, R.~Pemantle, and Y.~Peres.
\newblock A phase transition in random coin tossing.
\newblock {\em Ann.\ Probab.}, 29(4), 2001.

\bibitem[Ma99]{Ma}
H.~Matzinger. 
\newblock Reconstructing a three-color scenery by observing
it along a simple random walk path. 
\newblock {\em Rand.\ Struc.\ Algorithms}
15(2): 196--207, 1999. 

\bibitem[MR03a]{MR:Poly}
H.~Matzinger and S.~Rolles. 
\newblock Reconstructing a piece of scenery
with polynomially many observations.
\newblock {\em Stochastic Process.\ Appl.}
107(2): 289--300, 2003.

\bibitem[MR03b]{MR:Err}
H.~Matzinger and S.~Rolles. 
\newblock Reconstructing a random
scenery observed with random errors along a random walk path.
\newblock {\em Probab.\ Theory Related Fields}
125(4): 539--577,  2003.



\bibitem[Woe00]{W:RW}
W.~Woess.
\newblock {\em Random Walks on Infinite Graphs and Groups}, volume 138 of {\em
  Cambridge Tracts in Mathematics}.
\newblock Cambridge University Press, Cambridge, 2000.

\end{thebibliography}

\end{document}